\documentclass{article}
\pdfoutput=1
\usepackage{amsfonts, amsmath, amssymb, mathrsfs, verbatim} 
\title{$\beta$-Spaces: New Generalizations of Typically-Metric Properties \\ PREPRINT} 
\author{Annie Carter \\ Daniel Lithio \\ Robert Niichel \\ Tristan Tager} 
\newtheorem{THM}{Theorem}[section] 
\newtheorem{PROP}[THM]{Proposition}  
\newtheorem{LEMMA}[THM]{Lemma} 
\newtheorem{REMARK}[THM]{Remark} 
\newtheorem{NOTATION}[THM]{Notation} 
\newtheorem{DEFN}[THM]{Definition}  
\newtheorem{EXAMPLE}[THM]{Example} 
\newtheorem{COROLLARY}[THM]{Corollary} 
\newtheorem{CONJECTURE}[THM]{Conjecture} 
\newcommand{\thm}[1]{\begin{THM} #1 \end{THM} } 
\newcommand{\prop}[1]{\begin{PROP} #1 \end{PROP}} 
\newcommand{\lemma}[1]{\begin{LEMMA} #1 \end{LEMMA}} 
\newcommand{\remark}[1]{\begin{REMARK} #1 \end{REMARK}} 
\newcommand{\notation}[1]{\begin{NOTATION} #1 \end{NOTATION} } 
\newcommand{\defn}[1]{\begin{DEFN} #1 \end{DEFN}} 
\newcommand{\example}[1]{\begin{EXAMPLE} #1 \end{EXAMPLE}} 
\newcommand{\cor}[1]{\begin{COROLLARY} #1 \end{COROLLARY}} 
\newcommand{\conjecture}[1]{\begin{CONJECTURE} #1 \end{CONJECTURE}} 
\newcommand{\Thm}[2]{\begin{THM}[#1] #2 \end{THM}} 
 
\newcommand{\Lemma}[2]{\begin{LEMMA}[#1] #2 \end{LEMMA}}

\newcommand{\Defn}[2]{\begin{DEFN}[#1] #2 \end{DEFN} } 
\newcommand{\Example}[2]{\begin{EXAMPLE} [#1] #2 \end{EXAMPLE}}

\newcommand{\Pf}[1]{\emph{Proof:} 
\\ 
\par 
#1 $\quad \QED$ \newline \par} 
\newcommand{\QED}{\blacksquare} 

\renewcommand{\P}{\wp} 
 
\newcommand{\N}{\mathbb{N}} 
 
\newcommand{\seqq}[1]{( #1_{i} )_{i=1}^{\infty}}

\newcommand{\Z}{\mathbb{Z}} 
\renewcommand{\L}{\mathcal{L}} 
\newcommand{\Rstar}{{}^{*}\mathbb{R}}

\newcommand{\R}{\mathbb{R}} 
\newcommand{\Q}{\mathbb{Q}}

\begin{document} 
\maketitle 
\abstract{It is well-known that point-set topology (without additional structure) lacks the capacity to generalize the analytic concepts of completeness, boundedness, and other typically-metric properties.  The ability of metric spaces to capture this information is tied to the fact that the topology is generated by open balls whose radii can be compared.  In this paper, we construct spaces that generalize this property, called $\beta$-spaces, and show that they provide a framework for natural definitions of the above concepts.  We show that $\beta$-spaces are strictly more general than uniform spaces, a common generalization of metric spaces.  We then conclude by proving generalizations of several typically-metric theorems, culminating in a broader statement of the Contraction Mapping Theorem.} 
\newpage 
\section{Introduction} 
\par 
A $\beta$-structure on a topological set $X$ is a basis for the topology that consists of balls of varying radii.  We take a set $R$ of radius values, and a function $\beta$ that takes $x \in X$ and $r \in R$ to (intuitively) the ball of radius $r$ around the point $x$.  This allows us to describe neighborhoods around different points that are \emph{equivalently small}. It is this idea that allows for a generalization of a sequence being Cauchy, among other properties. 
\notation{We denote the power set of $S$ by $\P(S)$.} 
\Defn{$\beta$-Space}{\label{def:beta}A Topological space $X$, together with a set $R$ and a function $\beta: X \times R \to \P(X)$, is called a \emph{$\beta$-Space} if for all $r, s \in R$ and $x, y \in X$, 
\begin{enumerate} 
\item $x \in \beta(x, r)$ 
\item $\beta(x, r)$ is open 
\item For any open set $U$ and any point $z \in U$, there exists a $t \in R$ such that $\beta(z, t) \subseteq U$ 
\item For any $r \in R$ there exists an $s \in R$ such that, for all $x, y \in X$, if $x \in \beta(y, s)$, then $\beta(y, s) \subseteq \beta(x, r)$. 
\end{enumerate}} 
We call $R$ the radial values and $\beta$ the \emph{$\beta$-function}. Given $r \in R$, a radial value $s$ satisfying condition (4) above is called a \emph{swing value of $r$}, and the set $\beta(y, s)$ a \emph{swing set}.  Swing values and swing sets are an essential component of $\beta$-spaces.  Given $r \in R$, a swing value $s$ is essentially a value \emph{at most half of $r$}.  As we will see below, this intuition is captured exactly in metric spaces, which have a natural $\beta$-structure. \\~~\\ 
\par 
Although here we define $\beta$-spaces as having a topology -- for which, by condition (3) above, the sets $\beta(x, r)$ form a basis -- we can equivalently consider that $X$ does not begin with a topology, but rather has a natural topology induced by taking $\{ \beta(x, r) \; : \; x \in X, r \in R \}$ as a basis.  \\~~\\ 
\par

\Defn{Pseudo Field-Metric}{A \emph{pseudo field-metric space} is a linearly-ordered field $F$, a set $X$, and a function $d:X \times X \to F$ satisfying \begin{enumerate} 
\item $d(x, y) = d(y, x)$ 
\item $d(x, z) \leq d(x, y) + d(y, z) $ 
\end{enumerate}  
If in addition $d$ satisfies $d(x, y) = 0$ if and only if $x=y$, we drop the ``pseudo'' and call $d$ a \emph{field-metric}.  We often write a field-metric space as an ordered triple $(X, F, d)$.} 
We note that as in the case of metric spaces, since $d$ is symmetric and satisfies the triangle inequality, it must take values in $F^{\geq}$. 
\Thm{Field-Metrics Have a Natural $\beta$-Structure}{Let $(X, F, d)$ be a pseudo field-metric space. Then $(X, F, d)$ has a natural $\beta$-structure, and if in addition $(X, F, d)$ is a field-metric, the topology induced by the $\beta$-structure is Hausdorff.} 
\Pf{First, we topologize $(X, F, d)$ using the basis \[U = \left\{ B(x, r) \; : \; x \in X, r \in F^+ \right\} \] where, as usual, we define $B(x, r) = \{ y \in X \; : \; d(x, y) < r \}$. Now we define our radial set $R = F^+$ to be all positive elements in $F$. Finally, we define our $\beta$-function $\beta(x, r) = B(x, r)$.  We now need to show that these definitions meet all the criteria of a $\beta$-space as given in the above definition. 
\begin{enumerate} 
\item Obviously for any $x \in X$ and for any $r \in R$, since $r$ is necessarily positive, $x \in B(x, r) = \beta(x, r)$. 
\item Clearly $\beta(x, r)$ is open, since $B(x, r)$ is always open by definition. 
\item Since the basis for the topology is the open balls, then for any open set $U$ and any $x \in U$, there exists some $r \in R$ such that $B(x, r) \subseteq U$, which implies that $\beta(x, r) \subseteq U$. 
\item Given $r \in R$, to find a swing value for $r$ we simply choose $s = r/2$. We note that $r/2$ is well-defined, since $F$ is necessarily a characteristic 0 field, so that $2 \in F$, and $F$ closed under nonzero-division implies $r/2 \in F$. Now, let $x$ and $y$ be such that $x \in \beta(y, s) = B(y, r/2)$. Then for any $z \in \beta(y, s)$, $z \in B(y, r/2)$ implies that $d(y, z) < r/2$. Now since $x \in B(y, r/2)$ we know that $d(x, y) < r/2$, which by the triangle inequality guarantees that $d(x,z) < r/2 + r/2 = r$, so that $z \in B(x, r)$. But since this is true for all $z$, we know that $B(y, r/2) \subseteq B(x, r)$, or rather that $\beta(y, s) \subseteq \beta(x, r)$ as desired. 
\end{enumerate} 
Finally, we note that $X$ is Hausdorff in the induced topology when $d$ satisfies that $d(x, y) = 0$ if and only if $x=y$.  To see this, choose any $x, y \in X$ with $x \neq y$.  Then necessarily $d(x, y) \neq 0$, so that $r = d(x, y) / 2$ must also be nonzero.  Then $\beta(x, r) \cap \beta(y, r) = \varnothing$ are two disjoint open neighborhoods of $x$ and $y$ respectively.} 
  
\Example{The HyperReals}{The hyperreals provide an immediate example of $\beta$-space whose $\beta$-structure is induced by a field-metric. First, note that $\Rstar$ is a linearly-ordered field, where the linear ordering is given by our choice of ultrafilter. We then naturally define the field metric $d : \Rstar \times \Rstar \to \Rstar$ by \[d(x, y) \; = \; | x - y | \; = \; \left\{ \begin{array}{cc} x - y & x \geq y \\ \\ y - x & y > x \end{array} \right. \]   It is trivial to see that the function $d$ satisfies all the restrictions of a field-metric.} 
  
$\beta$-spaces induced by a field-metric enjoy many properties not guaranteed to hold for arbitrary $\beta$-spaces.  One of these is \emph{symmetry}, which is defined in Section 2 below.  Although $\beta$-spaces don't always enjoy this property, condition (4) of $\beta$-spaces does guarantee a weaker form of symmetry, whereby if $r \in R$ and $s$ is a swing value of $r$, then if $y \in \beta(x, s)$, necessarily $x \in \beta(y, r)$.  \\~~\\ 
\par 
  
\prop{Let $X$ be a $\beta$-space, and let $\seqq{a} \subseteq X$ be a sequence.  Then $a_i $ converges to $a$ if and only if for all $r \in R$ there exists some $N \in \N$ such that, for all $n \geq N$, $a_n \in \beta(a, r)$.} 
\Pf{$\Rightarrow$: 
Let $r\in R$ be given. By parts 1 and 2 of Definition \ref{def:beta}, $\beta(a,r)$ is a neighborhood of $a$. Therefore, since $a_i$ converges to $a$ there exists  some $N \in \N$ such that, for all $n \geq N$, $a_n \in \beta(a, r)$. 
\par 
$\Leftarrow$: Let a neighborhood $U$ of $a$ be given. By part 3 of Definition \ref{def:beta}, there is some $t\in R$ such that $\beta(a,t)\subset U$. By hypothesis, there is some $N \in \N$ such that, for all $n \geq N$, $a_n \in \beta(a, t) \subset U$. Therefore, $a_i$ converges to $a$. }

\Defn{Cauchy}{Let $X$ be a $\beta$-space, and let $\seqq{a} \subseteq X$. We say that $a_{i}$ is a \emph{Cauchy sequence}, or is \emph{Cauchy}, if for all $r\in R$, there exists an integer $N$ such that, for all $n, m \geq N$, $a_{m} \in \beta\left(a_n,r\right)$.} 
\Defn{Complete}{Let $X$ be a $\beta$-space. We say that $X$ is \emph{complete} if every Cauchy sequence in $X$ converges to a point in $X$.} 
\thm{The formal Laurent series $L$ in one variable over any linearly ordered field $F$ forms a complete first-countable $\beta$-space.} 
\Pf{It is well-known that $L$ is a field.  We now linearly order $L$ via the linear order on $F$ by taking \[ \sum_{i=-n}^{\infty} a_i x^i > 0 \quad \Leftrightarrow \quad a_{-n} > 0 \]  This establishes that $L$ is a field-metric.  To see that $L$ is first-countable in the induced topology, we notice that for any $f \in L$ with $f > 0$, there exists an $n \in \N$ such that $0 < x^n < f$.  Therefore for each $y \in L$, the set $\{ \beta(y, x^n) \; : \; n \in \N\}$ forms a countable collection of neighbourhoods with the required property.  \\~~\\ 
\par 
We now show that $L$ is complete in the induced $\beta$-.  Let $f_i(x) \in L$ be Cauchy.  For any $k \in \Z$, let $N_k$ be such that, for all $n, m \geq N_k$, $f_n \in \beta(f_m, x^k)$.  We can write \[ f_{n}(x) = \sum_{i=-\infty}^{\infty} a^n_i x^i \] where for each $n$ all the $a^n_i$ terms are 0 for all $i$ sufficiently small.  We now claim that for all $n, m \geq N_k$, $a^n_i = a^m_i$ for all $i < k$.  Suppose not.  Then there is an $n_0, m_0 \geq N_k$ and a smallest $j < k$ where $a^{n_0}_j > a^{m_0}_j$, without loss of generality.  Since given $n_0$ and $m_0$, $j$ was chosen to be the smallest index where the $a$'s differ, then we know that \[f_{n_0} - f_{m_0} = ( a^{n_0}_j - a^{m_0}_j ) x^j + \textrm{ other terms } > x^{j+1} \geq x^k  \]  which contradicts the fact that \[ f_n \in \beta(f_m, x^k) = \{ g \; : \; |g-f_m| < x^k \} \]  This necessarily establishes that there is an $\alpha \in \Z$ such that $a^n_i = 0$ for all $i< \alpha$ for any $n$.  We now define \[ g(x) = \sum_{i=\alpha}^{\infty}  a^{N_{i+1}}_i x^i \] and we claim that $f_i(x) \to g(x)$.  Pick any $r \in R$, and let $k \in \N$ be such that $x^k \leq r$.  Now, we know that for all $n, m \geq N_{k+1}$, $a^n_i = a^m_i$ for all $i < k$.  This gives us that, for all $n \geq N_{k+1}$, $a^n_i = a^{N_{k+1}}_i$ for all $i \leq k$.  This implies that we can write \begin{eqnarray*} f_n(x) & = & \sum_{i=\alpha}^{\infty} a^n_i x^i \\ & = & \sum_{i \leq k} a^n_i x^i \; + \; \sum_{i > k} a^n_i x^i \\ & = & \sum_{i \leq k} a^{N_{k+1}}_i x^i \; + \; \sum_{i > k} a^n_i x^i \end{eqnarray*}  Similarly, we can write \begin{eqnarray*} g(x) & = & \sum_{i \leq k} a^{N_{k+1}}_i x^i \; + \; \sum_{i>k} a^{N_{i+1}} x^i \end{eqnarray*}  Therefore, \begin{eqnarray*} |g - f| & = & \sum_{i > k} \left| a^{N_{i+1}}_i - a^n_i \right| x^i \\ & < & x^k \end{eqnarray*}  Therefore, for all $n \geq N_{k+1}$, $ f_n(x) \in \beta(g(x), x^k) \subseteq \beta(g(x), r) $ and so we have that $f_n(x) \to g(x)$ as desired.  Since $f_n(x)$ was an arbitrary Cauchy sequence, it follows that $L$ is complete.} 
The most important consequence of the above theorem, aside from its future relevance to fractal construction, is that it shows the existence of complete spaces that are non-metrizable.  To see that $L$ is non-metrizable, we note that $L$ is a nonarchimedean field, where if $L$ is taken over the reals, the (infinitesimal) element $x \in L$ has the property that $0 < x < r$ for all $r \in \R^+$.

\Defn{Totally Bounded}{Let $X$ be a $\beta$-space, and let $A \subseteq X$. We say $A$ is \emph{totally bounded} if, for all $r\in R$, there exists a finite set of points $\{x_i\}_{i=1}^n$ such that $A\subseteq\left({\bigcup_i \beta(x_i,r)}\right)$.} 
\lemma{\label{lem:Cauchysub}Let $X$ be a $\beta$-space, and let $\seqq{a}$ be a Cauchy sequence in $X$. Then any subsequence of $a_{i}$ is also a Cauchy sequence.} 
  
\Pf{ Trivial.} 
\prop{\label{prop:IfConvThenCauchy}Let $X$ be a $\beta$-space, and let $\seqq{a}$ be a sequence in $X$. If $a_{i}$ converges, then $a_{i}$ is Cauchy.} 
\Pf{Let $r\in R$ be any radial value, and let $s\in R$ be swing value of $r$. Since $a_{i}$ converges to some point $a \in X$, then there exists an $N$ such that, for all $m,n \geq N$, $a_{n} \in \beta(a,s)$ and $a_m \in \beta(a,s)$. Thus, since $s$ is a swing value of $r$, $\beta(a,s) \subseteq \beta(a_n,r)$ for all $n \geq N$. However, we have already asserted that $a_m\in\beta(a,s)\subseteq\beta(a_n,r)$, and since $r$ was chosen arbitrarily, $\seqq{a}$ is Cauchy.} 
\prop{\label{prop:Cauchy}Let $X$ be a $\beta$-space. Then $\seqq{a}$ converges to a point $b \in X$ if and only if for all $r\in R$, there exists an $N$ such that, for all $n \geq N$, $b \in \beta(a_n,r)$.} 
\Pf{$\Rightarrow$: 
Suppose $\seqq{a}$ converges to $b \in X$. Let $r\in R$ be given, and let $s\in R$ be the swing value of $r$. Since $a_i$ converges to $b$, there exists an $N$ such that, for all $n \geq N$, $a_{n} \in \beta(b,s)$. But since $s$ is a swing value for $r$ and $a_{n} \in \beta(b,s)$, then $\beta(b,s)\subset \beta(a_n,r)$. Therefore, $b\in \beta(a_n,r)$. 
\par 
$\Leftarrow$: 
Let $r\in R$ be given and let $s$ be a swing value for $r$. By hypothesis, choose $N$ so that, for all $n \geq N$, $b \in \beta(a_n,s)$. Since $s$ is a swing value for $r$ and $b \in \beta(a_n,s)$, we have that $\beta(a_n,s) \subseteq \beta(b,r)$, and thus that $a_{n} \in \beta(b,r)$. Therefore, $a_n$ converges to $b$.} 
  
\prop{\label{prop:Cauchy2}Let $X$ be a $\beta$-space.  Suppose $\seqq{a}$ is a sequence with $a_i \to a$.  Let $r \in R$ and $s$ be a swing value of $r$.  If $N$ is such that, for all $n, m \geq N$, $a_n \in \beta(a_m, s)$, then for all $n \geq N$, \begin{enumerate} 
\item $a_n \in \beta(a, r)$ 
\item $a \in \beta(a_n, r)$ 
\end{enumerate}} 
\Pf{Let $M \geq N$ be such that for all $k \geq M$ and all $n \geq N$, $a, a_n \in \beta(a_m, s)$.  Then $a \in \beta(a_n, r)$ and $a_n \in \beta(a, r)$ as desired.} 
  
The above proposition is intuitively clear: we can alternately characterize convergence in a $\beta$-space by requiring that the limit point be in most $\beta$-d neighborhoods of terms of the sequence. This alternate characterization will be useful in simplifying several forthcoming proofs. 
  
\prop{\label{prop:Cauchyconverge}Let $X$ be a $\beta$-space. If $\seqq{a} \subseteq X$ is Cauchy and has a subsequence $(a_{n_i})$ that converges to some $a \in X$, then $a_{i}$ itself converges to $a$.} 
\Pf{ 
Let $r\in R$, and let $s$ be a swing value for $r$. Since $\seqq{a}$ is Cauchy, there exists an $N_1$ such that if $n,m\geq N_1$, $a_m\in\beta(a_n,s)$. Since the subsequence $a_{n_i}$ converges to $a$, by Proposition \ref{prop:Cauchy} there exists an $N_2$ such that if $n_i\geq N_2$, $a\in\beta(a_{n_i},s)$. Let $N=\max\{N_1,N_2\}$. Then if $n,n_i\geq N$, we have $a_n,b\in\beta(a_{n_i},s)$. Since $s$ is a swing value for $r$, we have $\beta(a_{n_i},s)\subseteq\beta(a_n,r)$. Thus $a\in\beta(a_n,r)$ for all $n\geq N$, so by Proposition \ref{prop:Cauchy}, $a_i$ converges to $a$. 
}

For any $\beta$-space $(X, R, \beta)$, we note that we can define a partial order $\leq$ on $R$ by $r \leq s$ if and only if \[ \beta(x, r) \subseteq \beta(x, s) \] for all $x \in X$.  Many upcoming definitions and theorems will make use of this partial order. 
\Defn{Ordered 
}{A $\beta$-space $(X, R, \beta)$ is \emph{ordered} if the partial order $\leq$ defined on $R$ is a linear order.} 

\prop{For any field-metric space, the induced $\beta$-structure is ordered
.} 
\Pf{If $(X, F, d)$ is our field metric space, then in the induced $\beta$-structure we take $R = F^+$, which, being a subset of $F$, is linearly ordered by hypothesis. Now suppose that $r, s \in R$ and $r \leq s$. Then for any $x$, \[ \beta(x, r) = B(x, r) \subseteq B(x, s) = \beta(x, s) \] as desired.}

The following lemma and theorem closely follow the proofs for the equivalent results for uniform spaces provided in "General Topology" by Willard. Refer to that book for the definitions of net, subnet, and Cauchy net. Nets must be used so that these results can be applied to spaces that are not first countable. First, we give an equivalent characterization of totally bounded. 
\prop{\label{prop:totallyboundedalternate}A $\beta$-space $X$ is totally bounded if and only if for every $r\in R$ there exists a finite covering of $X$, $\{ U_i \}$ such that for every $i$, $x\in U_i$ implies that $U_i\subset \beta (x,r)$.} 
\Pf{For the forwards direction, let $r\in R$ be given. Since $X$ is a totally bounded $\beta$-space, there exists a swing value $s$ for $r$ and a finite cover consisting of balls of radius $s$. This cover meets the necessary conditions for the $U_i$'s. \\ 
The backwards direction is trivial.} 
\lemma{\label{lemma:Cauchytotallybdd} A $\beta$-space $X$ is totally bounded if and only if each net in $X$ has a Cauchy subnet.} 
\Pf{For the forwards direction, let $(x_\lambda)$ be a net in the totally bounded space $X$. Since $X$ is totally bounded, for any $r\in R$ thereat least onet $U_r \in X$ (as in Proposition \ref{prop:totallyboundedalternate}) such that $(x_\lambda)$ is frequently in $U_r$. Now, define $\Gamma = \{ (\lambda , r) ~|~ r\in R ~\text{and}~ x_\lambda \in U_r \}$, directed by $(\lambda _1 , r_1) \leq (\lambda _2 , r_2)$ if and only if $\lambda _1 \leq \lambda _2$ and $U_{r_1} \supset U_{r_2}$. Now, for each $(\lambda ,r)\in \Gamma$ define $x_{(\lambda ,r )}=x_\lambda$. Then, $(x_{(\lambda ,r )})$ is a subnet of $(x_\lambda)$. We show that $(x_{(\lambda ,r )})$ is Cauchy. Given $r_0\in R$, choose $\lambda _0 \in \Lambda$ so that $(\lambda _0 , r_0)\in \Gamma$. Then 
\[ (\lambda ,r ) , (\lambda ' , r') \geq (\lambda _0 , r_0) \Rightarrow x_\lambda ,x_{\lambda '} \subset U_{r_0} \Rightarrow x_\lambda \subset \beta (x_{\lambda '} , r_0 ) \] 
so that $(x_{(\lambda ,r )})$ is a Cauchy subnet of $(x_\lambda)$. \\ 
For the backwards direction, suppose that $X$ is not totally bounded. Then there exists $r\in R$ such that no finite cover ${U_i}$ as in Proposition \ref{prop:totallyboundedalternate} exists. Then, take $s$ to be a swing value of $r$. Now, since any cover $\{ \beta (x_i,s) \}$ is such a cover of $X$, there cannot be a finite sequence $\{x_i\}$ such that $\{\beta (x_i , s) \}$ is a cover of $X$. So, by induction, we can construct a sequence $x_1,x_2,\ldots$ such that $x_n \not \in \beta (x_i , r)$ for any $i <n$. Then, $(x_i)$ can not have a Cauchy subnet. 
} 
\thm{A $\beta$-space is compact if and only if it is complete and totally bounded} 
\Pf{For the forwards direction, let $X$ be compact and $(x_\lambda)$ a Cauchy net in $X$. Since any net in a compact space $X$ has a cluster point, we know that $(x_\lambda)$ has some cluster point $x$. Since $(x_\lambda)$ is Cauchy, it must converge to $x$. Therefore, $X$ is complete. Since $X$ is compact, every net has a convergent, and hence Cauchy, subnet, we know that $X$ is also totally bounded. \\ 
For the backwards direction, consider that by Lemma \ref{Cauchytotallybdd}, every net has a Cauchy subnet. Then, since $X$ is complete, this subnet is convergent. Therefore, $X$ is compact. } 

Although it is encouraging to see this theorem proven for $\beta$-spaces, this theorem will be overly restrictive for many future uses. For example, if we look at the natural $\beta$-structure on the hyperreals, we can see that there are very few compact sets of the sort which we may intuitively expect. In particular, no interval $[a, b] \subseteq \Rstar$ will be compact or totally bounded, since for any such interval, there is some $x \in \Rstar$ that is infinitesimal compared to $(b-a)$; and thus, there can be no finite cover of that interval by balls of radius $x$.

\section{Relationship to Uniform Spaces} 
Since uniform spaces are the most commonly used generalizations of metric spaces, it is of obvious interest to investigate the relationship between $\beta$-spaces and uniform spaces.  It turns out that it is fairly trivial to characterize the relationship precisely.  We begin with the following definitions. 
\Defn{Symmetric}{A $\beta$-space $(X, R, \beta)$ is called \emph{symmetric} if, for all $x, y \in X$ and all $r \in R$, $y \in \beta(x, r)$ if and only if $x \in \beta(y, r)$.} 
\defn{If a $\beta$-space $(X, R, \beta)$ is said to be \emph{closed under intersections} if it has the property that for all $r, s \in R$ there exists a $t \in R$ where \[ \beta(x, r) \cap \beta(x, s) = \beta(x, t) \] for all $x \in X$.} 
\thm{There is a natural one-to-one correspondence between symmetric $\beta$-spaces that are closed under intersections, and uniform spaces.} 
\Pf{First we establish that symmetric $\beta$-spaces that are closed under intersections induce a uniformity.\\~~\\ 
\par 
Given our $\beta$-space $(X, R, \beta)$, construct $\Phi$ as follows. For any $r \in R$, define \[ U_r = \{ (x, y) \; : \; y \in \beta(x, r) \} \] Now let $\Phi'$ be the collection of all $U_r$ for all such radial values $r$. Finally, let $\Phi$ be the set of all $V \subseteq X \times X$ such that there is a $U \in \Phi'$ with $U \subseteq V$. We now need to show that $\Phi$ satisfies all the conditions on uniform spaces.\\~~\\ 
\par 
\begin{enumerate} 
\item For any $U \in \Phi$ there exists an $r \in R$ such that $U_r \subseteq U$.  Since $x \in \beta(x, r)$, we know that $(x, x) \in U_r$ for every $x$, and so we know that $\Delta X \subseteq U_r \subseteq U$ as desired. 
\item Given any $U \in \Phi$, and a $V \subseteq X \times X$ such that $U \subseteq V$, we want to know that $V \in \Phi$. By construction of $\Phi'$ there must exist an $r \in R$ such that $U_r \subseteq U$, and therefore $U_r \subseteq V$. Then by construction of $\Phi$, $V \in \Phi$ as desired. 
\item Let $U, V \in \Phi$.  Then there exist $r, s \in R$ such that $U_r \subseteq U$ and $U_s \subseteq V$.  Since our $\beta$-space is closed under intersections, there exists a $t \in R$ such that $\beta(x, r) \cap \beta(x, s) = \beta(x, t)$ for all $x \in X$.  We claim that $U_r \cap U_s = U_t$.  To see this, we simply note that \[ U_r \cap U_s = \{ (x, y) \; : \; y \in \beta(x, r) \cap \beta(x, s) \} = \{ (x, y) \; : \; y \in \beta(x, t) \} \]  But since $U_t = U_r \cap U_s$, necessarily $U_t \subseteq U \cap V$.  Since by definition of $\Phi$, $U_t \in \Phi$, the fact that $\Phi$ is closed under supersets gives us that $U \cap V \in \Phi$. 
\item Given any $U \in \Phi$, we want to show that there exists a $V \in \Phi$ such that, whenever $(x, y), (y, z) \in V$, then $(x, z) \in U$. Fixing our $U \in \Phi$, we know there exists an $r \in R$ such that $U_r \subseteq U$. Let $s \in R$ be a swing-value of $r$. Now set $V = U_s$. Suppose that $(x, y), (y, z) \in V$. Since $(x, y) \in V$, that implies that $y \in \beta(x, s)$.  Since our $\beta$-space is symmetric by hypothesis, this implies that $x \in \beta(y, s)$. Since $s$ is a swing value of $r$, this gives us that $\beta(y, s) \subseteq \beta(x, r)$. Now since $(y, z) \in V$, we know that $z \in \beta(y, s)$, so that $z \in \beta(x, r)$ necessarily. Thus, $(x, z) \in U_r \subseteq U$ as desired. 
\item Fix any $U \in \Phi$. Then we claim that $U^{-1} = \{ (y, x) \; : \; (x, y) \in U \} \in \Phi$ as well. To see this, let $r$ be such that $U_r \subseteq U$. Then $(x, y) \in U_r$ iff $y \in \beta(x, r)$ iff $x \in \beta(y, r)$ iff $(y, x) \in U_r$, so that $U_r = U_r^{-1}$. But then $U_r \subseteq U^{-1}$, so that, since $\Phi$ is closed under supersets, $U^{-1} \in \Phi$ as desired. 
\end{enumerate} 
\par 
Now we show that any uniform space induces a symmetric $\beta$-space closed under intersections.\\~~\\ 
\par 
Given our set $X$ and our uniformity $\Phi$, let \[ R = \{ U \cap U^{-1} \; : \; U \in \Phi \} \]   where $U^{-1} = \{ (y, x) \; : \; (x, y) \in U \}$ as usual.  Since the uniformity $\Phi$ is closed under these ``inverses'', and also under intersections, $U \in \Phi$ implies that $U \cap U^{-1} \in \Phi$, so for the remainder of the proof we will simply refer to an element $U \in R$, recalling when necessary that all such $U$ have the property that $(x, y) \in U$ if and only if $(y, x) \in U$.  Now for any $x \in X$ and $U \in R$, we define \[ \beta(x, U) = U[x] = \{y \; : \; (x, y) \in U \}\]  Noting that the topology on $X$ is necessarily the topology induced by the uniformity, we now need to show that these definitions satisfy all the axioms of $\beta$-spaces.\\~~\\ 
\par 
\begin{enumerate} 
\item The fact that $x \in \beta(x, U)$ follows automatically, since by definition of uniform spaces, $\Delta X \subseteq U$ for all $U \in \Phi$. Thus, for any $x \in X$, $(x, x) \in U$, so that $x \in U[x] = \beta(x, U)$ as desired. 
\item The fact that $\beta(x, U)$ is open is also automatic, since the topology on $X$ is generated by all the sets $U[x] = \beta(x, U)$.  
\item We now want to show that the $\beta$-sets form a basis -- specifically, that for any open set $W$ and any $x \in W$, there exists a $U \in R$ such that $\beta(x, U) \subseteq W$. This is also automatic; it follows immediately from the way that the Uniform structure induces a topology. Specifically, the Uniform topology declares that a set $W$ is open if for any $x \in W$ there exists a $U \in \Phi = R$ such that $\beta(x, U) = U[x] \subseteq W$, which is precisely what we want. 
\item Finally, we want to show that we have our swing values -- that for any $U \in R$, there exists a $V \in R$ such that, if $x \in \beta(y, V)$, then $\beta(y, V) \subseteq \beta(x, U)$. Now having fixed $U$, pick $V' \in \Phi$ to be such that if $(x, y), (y, z) \in V'$, then $(x, z) \in U$, as guaranteed by condition 4 of Uniform Spaces. Then let $V = V' \cap (V')^{-1} \in R$.  Now suppose that $x \in \beta(y, V) = V[y]$. This implies that $(y, x) \in V$ which implies that $(x, y) \in V$. We want to show that $\beta(y, V) \subseteq \beta(x, U)$. Let $z \in \beta(y, V)$. Then $(y, z) \in V$. But now we have that $(x, y), (y, z) \in V$, which implies that $(x, y), (y, z) \in V'$, which implies that $(x, z) \in U$. Thus we have that $z \in \beta(x, U) = U[x]$, which completes the proof. 
\end{enumerate} 
Having shown that the uniformity induces a $\beta$-structure, we now need to show that it in fact induces a symmetric $\beta$-space closed under intersections.  We show symmetry first.  Let $x, y \in X$ and $U \in R$, and suppose that $y \in \beta(x, U)$.  Then $(x, y) \in U$.  But since $U \in R$, we know that this implies $(y, x) \in U$, so that $x \in U[y] = \beta(y, U)$ as desired.  \\~~\\ 
\par 
Finally, to see that our $\beta$-space is closed under intersections, let $U, V \in R$.  First we claim that $Z = U \cap V \in R$.  To see this, it suffices to show that $(x, y) \in Z$ if and only if $(y, x) \in Z$.  But this is trivial, since if $(x, y) \in Z$, then $(x, y) \in U$ and $(x, y) \in V$.  Since $U, V \in R$, this implies that $(y, x) \in U$ and $(y, x) \in V$, so that $(y, x) \in U \cap V = Z$ as desired.  Now we simply note that \[ \beta(x, U) \cap \beta(x, V) = U[x] \cap V[x] = Z[x] = \beta(x, Z) \] which completes the proof.} 
It should be noted that, in the above proof, the condition of uniform spaces that whenever $U \in \Phi$ and $V \subseteq X \times X$ with $U \subseteq V$, then $V \in \Phi$, was completely unnecessary in showing that a uniform space induced a symmetric $\beta$-space closed under intersection.  Furthermore, showing that a symmetric $\beta$-space closed under intersection induces a uniform space only required a brief nod to this condition.  It is apparent from this that, were this condition removed completely from the definition of uniform spaces, the above proof would still hold.  \\~~\\ 
\par 
The above theorem provides the interesting (and in many ways unexpected) result that we can consider the category of uniform spaces to be properly included in the category of $\beta$-spaces.

\prop{A $\beta$-space $(X,R,\beta)$ is regular.} 
\Pf{Suppose that $(X,R,\beta)$ is not regular. Then there exists a point $x$ and a closed set $K$ such that for all open $U_x$ and $U_K$ with $x\in U_x$ and $K\subset U_K$ we have $U_x\cap U_K \neq \emptyset$. We now find a $U_x$ and $U_K$ to form a contradiction. \\
To pick a $U_x$, first note that $(X\backslash K)$ is open. Then, there exists an $r\in R$ such that $\beta (x, r) \subset (X\backslash K)$. We let $U_x = \beta (x,r)$. \\
Next, we choose $U_K$. To construct this set, first let $s$ be a swing value of $r$, and $t$ a swing value of $s$. We define
\[ U_K = \bigcup _{y\in K} \beta (y,t) .\]
Now, by our supposition, there exists $z\in X$ such that $z\in U_x \cap U_K$. Moreover, there exists $y_0\in K$ such that $z\in \beta (y_0,t)$. Then, since $t$ is a swing value of $s$, we know that $\beta(y_0 ,t) , \beta(x,t)\subset \beta(z,s)$. Next, since $s$ is a swing value of $r$, we have $\beta(z,s)\subset \beta(x,r)$. However, this gives us that $y_0 \in \beta(x,r)$, which is a contradiction with how we chose the value of $r$. Therefore, $(X,R,\beta)$ must be regular.
}

It is well-known that uniform spaces are precisely the completely regular topological spaces.  It is natural to ask whether the same is true of $\beta$-spaces, or whether there exist regular, non-completely-regular spaces that admit a $\beta$-structure compatible with their topology.  Due to the scarcity and complexity of such spaces, this has proven a hard question to answer, and remains an open question.  We therefore formulate the following conjecture:

\conjecture{There exist non-completely-regular $\beta$-spaces.}
 
\section{Additional Mechanics} 
We now wish to apply our notions of $\beta$-spaces to generalize the notion of fractals. Much of the machinery required to properly construct fractals in a Metric space setting falls naturally out of notions of compact sets and continuous functions. However, in a more general setting, we will find that the condition of compactness is often too restrictive to generate interesting or intuitive fractals -- and if we are no longer looking at compact sets, we are no longer guaranteed that the continuous image of our sets will have the same properties. Because of this, our construction will require some new tools, which we will carefully construct before we begin to look at more traditional fractal machinery. 
\par 
  
\Defn{Swing Sequence}{Given $r \in R$, a \emph{swing sequence} for $r$ is a sequence $\seqq{r} \subseteq R$ such that $r_1 = r$ and, for each $k$, $r_{k+1}$ is a swing value for $r_k$.} 
\Lemma{Geometric Series}{Let $(X, R, \beta)$ be a $\beta$-space.  For any $r \in R$, let $s$ be a swing value of $r$, and let $t$ be a swing value of $s$.  Further, let $\seqq{t}$ be a swing sequence for $t$.  Let $\seqq{a} \subseteq X$ be a sequence such that $a_{n+1} \in \beta(a_n, t_n)$ for all $n$.  Then for all $n$, $\beta(a_n, t_n) \subseteq \beta(a_1, r)$.} 
\Pf{Trivially the proposition holds for $n=1$.  The case where $n=2$ is a trivial argument using swing values.  Thus we will assume $n \geq 3$.  For convenience, we will use $t_{0} = s$.  Now since $t_n$ is a swing value of $t_{n-1}$ and $t_{n-1}$ is a swing value of $t_{n-2}$, necessarily $t_n$ is a swing value of $t_{n-2}$.  Thus since $a_n \in \beta(a_n, t_{n})$ and $a_n \in \beta(a_{n-1}, t_{n-1})$, necessarily \[ \beta(a_n, t_n), \beta(a_{n-1}, t_{n-1}) \subseteq \beta(a_n, t_{n-2}) \]  Now similarly, since $a_{n-1} \in \beta(a_{n-2}, t_{n-2})$ by hypothesis and $a_{n-1} \in \beta(a_n, t_{n-2})$ as established above, then since $t_{n-2}$ is a swing value of $t_{n-3}$, necessarily \[ \beta(a_{n-2}, t_{n-2}), \beta(a_n, t_{n-2}) \subseteq \beta(a_{n-1}, t_{n-3}) \]  Inductively suppose that we have \[ \beta(a_{n-m+1}, t_{n-m+1}), \beta( a_{n-m+3}, t_{n-m+1}) \subseteq \beta(a_{n-m+2}, t_{n-m}) \] for $m < n$.  This gives us that $a_{n-m+1} \in \beta(a_{n-m+2}, t_{n-m})$.  By hypothesis we have that $a_{n-m+1} \in \beta(a_{n-m}, t_{n-m})$.  Thus, since $t_{n-m}$ is a swing value of $t_{n-m-1}$, \[ \beta(a_{n-m}, t_{n-m}), \beta(a_{n-m+2}, t_{n-m}) \subseteq \beta(a_{n-m+1}, t_{n-m-1}) \]  Since this applies for all $m<n$, by induction, this gives us that \[ \beta(a_n, t_n) \subseteq \beta(a_n, t_{n-2}) \subseteq \beta(a_{n-1}, t_{n-3}) \subseteq \ldots \subseteq \beta(a_3, t_1) \subseteq \beta(a_2, t_0) = \beta(a_2, s) \]  Finally, we note that since $s$ and $t$ are swing values of $r$, and $a_2 \in \beta(a_1, t_1)$, then \[ \beta(a_1, t_1), \beta(a_2, s) \subseteq \beta(a_1, r) \] which completes the proof. } 
\cor{Let $(X, R, \beta)$ be a $\beta$-space, and let $r, s, t$, $\seqq{t}$, and $\seqq{a}$ be as in the statement of the Geometric Series Lemma.  Suppose further that $r$ is a swing value of some $p \in R$.  Then for any $n, m \in \N$, $\beta(a_n, t_n) \subseteq \beta(a_m, p)$.} 
\Pf{By the Geometric Series Lemma, $\beta(a_n, t_n), \beta(a_m, t_m) \subseteq \beta(a_1, r)$ for all $n, m \in \N$.  Since $r$ is a swing value for $p$ and $a_m \in \beta(a_1, r)$, necessarily \[ \beta(a_n, t_n) \subseteq \beta(a_1, r) \subseteq \beta(a_m, p) \]  as desired.} 
  
\remark{GRAPHIC HERE to justify Geometric Prop.  Graphic is already created.}

\Defn{Level Set}{Given a point $x \in X$ and a $r \in R$, define the \emph{level set about x of radius r}, or the \emph{r-level set about x}, to be \begin{eqnarray*} L(x, r) & = & \{ y \in X \; : \; \textrm{there exists a swing sequence } \seqq{r} \\ & & \textrm{ for } r \textrm{ such that } y \in \beta(x, r_i) \textrm{ for all } i \} \end{eqnarray*}} 
\lemma{\label{lemma:swingseq}If $(X, R, \beta)$ is an ordered $\beta$-space and $\seqq{r}$ and $\seqq{s}$ are two swing sequences for $r \in R$, then the sequence $\seqq{t}$ given by $t_i = \max \{ r_i, s_i \}$ is a swing sequence for $r$.} 
\Pf{First, we note that $t_1 = \max \{r_1, s_1 \} = \max\{ r, r \} = r$ as required.  Without loss of generality, suppose that $t_{n+1} = r_{n+1}$.  Then if $t_{n} = r_{n}$, the condition that $t_{n+1}$ is a swing value for $t_n$ is trivially satisfied.  If instead we had $t_{n} = s_{n}$, we note that this implies that $r_n \leq s_n$ in the ordering on $R$.  Since $r_{n+1}$ is a swing value for $r_n$, we know that when $y \in \beta(x, r_{n+1})$, we have $\beta(x, r_{n+1}) \subseteq \beta(y, r_n)$.  But since $r_n \leq s_n$, $\beta(y, r_n) \subseteq \beta(y, s_n)$ for all $y$, so that we have \[ \beta(x, t_{n+1}) = \beta(x, r_{n+1}) \subseteq \beta(y, s_n) = \beta(y, t_n)\] as desired.} 
\prop{\label{prop:levelequality}If $(X, R, \beta)$ is an ordered $\beta$-space, then for any $y \in L(x, r)$, $L(x, r) = L(y, r)$.} 
\Pf{We first show that whenever $y \in L(x, r)$, then $x \in L(y, r)$.  If $y \in L(x, r)$, then there exists a swing sequence $\seqq{r}$ for $r$ such that $y \in \beta(x, r_{i+1})$ for all $i$.  Therefore \[ x \in \beta(x, r_{i+1}) \subseteq \beta(y, r_i) \] for all $i$ as desired. \\~~\\ 
\par 
Now let $z \in L(x, r)$.  Then there exists a swing sequence $\seqq{s}$ for $r$ such that $z \in \beta(x, s_{i+1})$.  Now by Lemma \ref{lemma:swingseq} we know that $\seqq{t}$ defined by $t_i = \max \{r_i, s_i \}$ is also a swing sequence for $r$, and clearly it has the property that $ y, z \in \beta(x, t_{i+1}) $ for all $i$.  But then we have that \[ z \in \beta(x, t_{i+1}) \subseteq \beta(y, t_i) \] for all $i$.   This gives us that $z \in L(y, r)$, and hence that $L(x, r) \subseteq L(y, r)$. To see the reverse inclusion, we let $w \in L(y, r)$, and we let $\seqq{p}$ be a swing sequence for $r$ such that $w \in \beta(y, p_i)$ for all $i$.  Then since as established above, $x \in \beta(y, r_i)$, the symmetric argument suffices to establish the reverse inclusion.} 
  
Level sets are a very interesting and useful part of analysis on $\beta$-spaces.  The intuition behind the set $L(x, r)$ is that it is the set of all points infinitesimally close to $x$ relative to $r$.  In the pseudometric approach to uniform spaces, $L(x, r)$ corresponds to the set of all $y$ such that $\rho _\alpha(x, y) = 0$ for some $\rho _\alpha$ in the family of pseudometrics that defines the uniform space. 
  
\Defn{Level Equivalent}{Two radial values $r, s \in R$ are said to be \emph{level equivalent}, written $r =_L s$, if $L(x, r) = L(x, s)$ for all $x \in X$.} 
\Defn{Level Ordered}{A $\beta$-space $(X, R, \beta)$ is said to be \emph{level-ordered} if there exists a linear order $\leq_L$ on $R / =_L$ such that $s \leq_L r$ whenever $L(x, s) \subseteq L(x, r)$ for all $x \in X$.} 
\prop{An ordered $\beta$-space is level ordered.} 
\Pf{Take $r\leq s$ in $R$, the set of radial values.  We now want to show that $r\leq_L s$, i.e. that $L(x,r)\subseteq L(x,s)$ for all $x\in X$.  So, take $y\in L(x,r)$ for some $x$, and suppose that $\seqq{r}$ is a swing sequence for $r$ such that $y\in\beta(x,r_i)$ for all $i$.  If we now define $\seqq{r^\prime}$  by  $r_1^\prime=s$ and $r_j^\prime=r_j$ for all $j\geq2$, then we note that because the $\beta$-space is ordered, $\beta(x,r_2^\prime)\subseteq\beta(x,r)\subseteq\beta(x,s)$, and a similar argument also shows that $r_2^\prime$ is a swing value of $s$, and so $\seqq{r^\prime}$ is a swing sequence for $s$.  Clearly then, $y\in\beta(x,r_i^\prime)$ for all $i$, and $y\in L(x,s)$ as desired.} 
It is worth noting that although an ordered $\beta$-space is automatically level ordered, the converse does not hold.  In particular we can have a $\beta$-space where $X = L$ is the field-metric of formal Laurent series, and the ``balls'' are open balls of positive radius, together with open squares of positive side length.  It is easy to see that this $\beta$-space is not ordered, and yet it is level ordered. 
  
\Defn{Contraction}{Let $(X, R, \beta)$ be a $\beta$-space.  We say a function $f: X \to X$ is a \emph{contraction} if there exists an integer $N \in \N$ such that  
\begin{enumerate}  
\item $f\left( \beta(x, r) \right) \subseteq \beta\left( f(x), r) \right)$ for all $x \in X$ and $r \in R$ 
\item Each $r \in R$ has a swing value $s$ so that $f^N\left( \beta(x, r) \right) \subseteq \beta\left(f^N(x), s \right) $ 
\end{enumerate} 
We call $N$ the \emph{contraction degree} of $f$.} 
It is readily apparent that, if a $\beta$-space is induced by a field-metric, then the common notion of contraction implies the definition given above.  It is worth noting that a function $f: X \to X$ being a contraction is a restriction on the $\beta$-space as well as on the function $f$.  Intuitively, for any nonarchimedean linearly ordered field $F$, we would want $f(x) = x/2$ to be a contraction.  But in fact this will be the case if and only if $R$ contains ``enough'' values.  If $R = F^+$, then $f$ is a contraction -- but if, for example, $F$ is the field of formal Laurent series in one variable over $\R$ and we pick $R = \{ x^{\alpha} \; : \; \alpha \in \N \}$, then we get an equally valid $\beta$-structure where $f$ is no longer a contraction, simply because for any $r \in R$, the only swing values of $r$ are also infinitesimal relative to $r$, and $f$ cannot contract that severely in only finitely many iterates.

\Defn{\emph{r}-Cauchy}{Given an $r \in R$, a sequence $\seqq{a}$ is said to be \emph{r-Cauchy} if there exists a swing sequence $\seqq{r}$ of $r$ such that, for all $n \in \N$, there is an $M \in \N$ where $a_j \in \beta(a_k, r_n)$ for all $j, k \geq M$.} 
\Defn{\emph{r}-Converge}{Given an $r \in R$, a sequence $\seqq{a}$ is said to \emph{r-converge} to a point $x \in X$, written $a_i \xrightarrow{r} x$, if there exists a swing sequence $\seqq{r}$ for $r$ such that, for all $n \in \N$, there is an $M \in \N$ where $a_m \in \beta(x, r_n)$ for all $m \geq M$.  The sequence $\seqq{r}$ is called a \emph{swing sequence for the convergence to x}.} 
Informally the notion of $r$-convergence is the notion of convergence if we throw away all radial values infinitesimal relative to $r$.  This process can (and usually will) render a space non-Hausdorff even if it was Hausdorff originally, and so a sequence -- even an $r$-Cauchy sequence -- may $r$-converge to multiple points. 
\Defn{\emph{r}-Convergence Set}{Given a sequence $\seqq{a} \subseteq X$ and an $r \in R$, the \emph{r-Convergence Set} for $\seqq{a}$ is the set \[ A = \{ x \in X \; : \; a_i \xrightarrow{r} x \} \]} 
\example{Let $F$ denote the field of formal Laurent series in one variable over $\R$.  Let $a_i = 2^{-i}$.  Then although our sequence $\seqq{a}$ does not converge, it does $r$-converge for $r = 1$.  We establish the swing sequence $\seqq{r}$ by $r_i = 2^{-i+1}$.  Then $a_i \xrightarrow{r} x$ if and only if $x$ is infinitesimally close to 0.} 
  
We now state two propositions about $r$-Cauchy and $r$-convergent sequences that are the analog of propositions about Cauchy and convergent sequences stated in the beginning of this paper.  The proofs are omitted, as they are virtually identical. 
\prop{\label{prop:rConvCauchy}Let $X$ be a $\beta$-space, and let $\seqq{a}$ be a sequence in $X$. If $a_{i}$ $r$-converges, then $a_{i}$ is $r$-Cauchy.} 
\prop{\label{prop:rCauchy}Let $(X, R, \beta)$ be a $\beta$-space. Then $\seqq{a}$ $r$-converges to a point $b \in X$ if and only if for all $k \in \N$ there exists an $N$ such that, for all $n \geq N$, $b \in \beta(a_n, r_k)$.} 
\Defn{\emph{r}-Complete}{A $\beta$-space $(X, R, \beta)$ is called \emph{r-complete} if any $r$-Cauchy sequence $r$-converges to some $x \in X$.} 
\Defn{Radially Complete}{A $\beta$-space $(X, R, \beta)$ is called \emph{radially complete} if it is complete and $r$-complete for all $r \in R$.} 
It is trivial to see that a sequence $\seqq{a}$ is Cauchy if and only if it is $r$-Cauchy for all $r$, and that the sequence converges if and only if it $r$-converges for all $r \in R$.  Because of this, trivially a radially complete space is complete.  However, the converse need not hold.  One example is that the field of formal Laurent series over the rationals is complete but is not $1$-complete.  Another example that will be relevant in our discussion of the Contraction Mapping Theorem is as follows.  Let $\L$ be the field of formal Laurent series in one variable over $\R$, and let $X$ be the set of all points infinitesimally close to some $y \in (0, 1) \subseteq \R$.  We can view $X$ as the closed set $[0, 1] \subseteq F$ with the set of all $x$ infinitesimally close to either $0$ or $1$ removed.  Since the latter two sets are open, this is a closed set minus two open sets, and thus is closed.  Since $F$ is complete, necessarily $X$ is also complete.  However, $(X, R, \beta)$ is not a radially complete $\beta$-space, since the sequence $a_i = 1/n$ is $r$-Cauchy for $r=1$ -- that is, it is $1$-Cauchy -- but does not $1$-converge in $X$.

\Defn{Structured $\beta$-Space}{A $\beta$-space $(X, R, \beta)$ is called \emph{structured} if for any $x \in X$ and any finite collection $(y_i)_{i=1}^{n} \subseteq X$, there exists an $s \in R$ such that $(y_i)_{i=1}^n \subseteq \beta(x, s)$.  The space is called \emph{level structured} if it satisfies that if there is an $r \in R$ such that $(y_i)_{i=1}^n \subseteq L(x, r)$, and if there is a $y_k$ that can be separated from $x$ by disjoint open neighborhoods, then we can choose $s$ so that $s <_L r$.}

\example{$\L$ field of Laurent series, \[ X = \{ x \in \L \; : \; x \textrm{ is infinitesimally close to some  } y \in (0, 1) \subseteq \R \}\]  $f(x) = x/2$ is a contraction without a fixed point.  $X$ is complete but not radially complete.} 
\example{Let $X = \{x, y\}$ be a metric space with $d(x, y) = 1$.  Then if $f$ is the identity function, it is a contraction.  The space is not level-structured, so $f$ does not have a unique fixed point.} 
\example{Any linearly ordered field must be char 0, and so contains $\Q$ as a subfield.  Suppose a field-metric space $(X, F, d)$ is level structured in the induced $\beta$-.  Then any function $f : X \to X$ satisfying that there exists a $r \in [0, 1) \subseteq \Q \subseteq F$ with $d(f(x), f(y)) \leq r \cdot d(x, y)$ is a contraction.} 
\prop{\label{prop:contraction}Let $(X, R, \beta)$ be a structured $\beta$-space and let $f:X \to X$ be a contraction with contraction degree $N$.  For any $x \in X$, let $r$ be such that $f^1(x), f^2(x), \ldots, f^N(x) \in \beta(x, r)$.  Then the sequence $a_n = f^n(x)$ is $r$-Cauchy, and if $A$ is the $r$-convergence set for $(f^n(x))_{i=1}^{\infty}$, then $f(A) \subseteq A$.} 
\Pf{First we fix $x$ and show that $f^n(x)$ is $r$-Cauchy.  Given $r$ as above, set $r_1 = r$.  Inductively, if we have $r_n$ defined, let $r_{n+1}$ be a swing value of $r_n$ such that \[ f^N\left( \beta(y, r_n) \right) \subseteq \beta\left( f^N(y), r_{n+1} \right) \]  for all $y$.  Now we notice that since $x, f(x), \ldots, f^N(x) \in \beta(x, r_1)$, that necessarily $f^N(x), f^{N+1}(x), \ldots, f^{2N}(x) \in f^N\left( \beta(x, r_1) \right) \subseteq \beta(f^N(x), r_2)$, and in general we have \[ f^{jN+1}(x), f^{jN+2}(x), \ldots, f^{(j+1)N}(x) \in \beta(f^{jN}(x), r_{j+1})  \]  Fixing a $k \in \N$, the sequences $a_i = f^{(i+k+2)N}(x)$ and $t_i = r_{i+k+3}$ clearly satisfy the conditions of the corollary to the Geometric Series Lemma, where for the values $p, r, s, t$ respectively we take $r_{k+1}, r_{k+2}, r_{k+3}$, and $r_{k+4}$, noting that $t = t_1 = r_{k+4}$ as required.  The corollary then gives us that, for all $n, m \in \N$, \[ \beta\left(f^{(n+k+2)N}(x), r_{n+k+3}\right) \subseteq \beta\left(f^{(m+k+2)N}(x), r_{k+1}\right) \]   Now set $M = (k+3)N$.  For any $n \geq M$, we can write $n = Na + b$ where $0 \leq b < N$.  Since necessarily $a \geq k+3$, we have that  \[ f^n(x) = f^{Na+b} \in \beta\left( f^{Na}(x), r_{a+1} \right) \subseteq \beta\left( f^{(m+k+2)N}(x), r_{k+1} \right)  \]  But since for any $n \geq M$ we have $f^n(x) \in \beta\left( f^{(m+k+2)N}(x), r_{k+1} \right)$ and that $r_{k+1}$ is a swing value of $r_k$, necessarily \[ \beta\left( f^{(m+k+2)N}(x), r_{k+1}\right) \subseteq \beta\left( f^n(x), r_k \right) \] for all such $n$.  But this gives us that for all $n, m \geq M$, \[ f^m(x) \in \beta\left( f^n(x), r_k \right) \] which completes the proof that the sequence $f^n(x)$ is $r$-Cauchy. 
\\~~\\ 
\par 
Now we show that $f(A) \subseteq A$.  This is trivial if $A$ is empty, so suppose that $x^* \in A$.  Thus $f^n(x) \xrightarrow{r} x^*$, so by Proposition \ref{prop:rCauchy}, for all $k$ there exists an $N$ such that $x^* \in \beta\left(f^n(x), r_{k+2}\right)$ for all $n \geq N$.  Now since $r_{k+2}$ is a swing value of $r_{k+1}$, $x^* \in \beta\left( f^n(x), r_{k+2} \right)$ implies that \[x^* \in \beta\left(f^n(x), r_{k+2} \right) \subseteq \beta\left( x^*, r_{k+1} \right) \]  for all $n \geq N$.  But by definition of contractions, applying $f$ to the above statement, we have that \[ f\left(x^*\right) \in f\left( \beta\left( f^n(x), r_{k+2} \right) \right) \subseteq \beta\left( f^{n+1}(x), r_{k+2} \right) \subseteq \beta\left( x^*, r_{k+1} \right) \]  Now since $r_{k+1}$ is a swing value of $r_k$, this gives us that $ \beta(x^*, r_{k+1} ) \subseteq \beta\left(f\left(x^*\right), r_k \right)$, which finally implies that \[ f^{n+1}(x) \in \beta\left( f^{n+1}(x), r_{k+2} \right) \subseteq \beta(x^*, r_{k+1} ) \subseteq \beta\left(f\left(x^*\right), r_k \right) \]  so that, for all $n \geq N$, $f^{n+1}(x) \in \beta\left( f\left( x^* \right), r_k \right) $.  Thus, $f^n(x) \xrightarrow{r} f\left( x^* \right)$, so that $f \left( x^* \right) \in A$ as desired.} 

\Defn{Countably Level-Based}{A $\beta$-space $(X, R, \beta)$ is called \emph{countably level-based} if there is a collection $\seqq{r} \subseteq R$ such that, for all $s \in R$ there exists an $r_k \leq_L s$.  The collection $\seqq{r}$ is called a \emph{countable level base}.} 
\Defn{CDLB}{A $\beta$-space is called \emph{countably and discretely level-based}, or \emph{CDLB}, if any collection $\seqq{r} \subseteq R$ that satisfies $r_{i+1} <_L r_i$ is a countable level-base.}

\prop{\label{prop:levelequals}Let $(X, R, \beta)$ be an ordered $\beta$-space and let $\seqq{a}$ be a sequence with $r$-convergent set $A \subseteq X$.  Then $A = L(x, r)$ for all $x \in A$.} 
\Pf{First we show that for all $x \in A$, $L(x, r) \subseteq A$.  Let $y \in L(x, r)$.  Then there exists swing sequences $\seqq{r}, \seqq{s}$ for $r$ such that $y \in \beta(x, r_i)$ for all $i$, and such that for any $m \in \N$ there exists an $N_m \in \N$ where $a_i \in \beta(x, s_m)$ for all $i \geq N_m$.  Then by Lemma \ref{lemma:swingseq} there exists a swing sequence $\seqq{t}$ for $r$ such that $t_i = \max \{ r_i, s_i \}$, so that necessarily we have 
\begin{itemize} 
\item $y \in \beta(x, t_i)$ for all $i$ 
\item $a_i \in \beta(x, t_m)$ for all $i \geq N_m$ 
\end{itemize} 
Now, for all $i$, $y \in \beta(x, t_{i+1})$ implies that $\beta(x, t_{i+1}) \subseteq \beta(y, t_i)$.  Thus, for any $m$, \[ a_i \in \beta(x, t_{m+1}) \subseteq \beta(y, t_m) \] for all $i \geq N_{m+1}$.  Thus we have $a_i \xrightarrow{r} y$ and so $y \in A$ as desired.\\~~\\ 
\par 
Now we show that $A \subseteq L(x, r)$.  Let $y \in A$, and by Lemma \ref{lemma:swingseq}, let $\seqq{t}$ be such that for each $m \in \N$ there is a $N_m \in \N$ such that \[ a_i \in \beta(x, t_{m+2}) \cap \beta(y, t_{m+2}) \] for all $i \geq N_{m+2}$, and so that in particular $a_{N_{m+2}}$ is contained in the intersection.  But this gives us that \[ \beta(x, t_{m+2}), \beta(y, t_{m+2}) \subseteq \beta(a_{N_{m+2}}, t_{m+1}) \subseteq \beta(x, t_m) \] so that $y \in \beta(x, t_m)$ for all $m$, and so $y \in L(x, r)$ as desired.} 
\lemma{If $(X, R, \beta)$ is an ordered $\beta$-space, then if $s <_L r$, there exists a swing sequence $\seqq{r}$ for $r$ such that $s < r_n$ for all $n$.} 
\Pf{If $s <_L r$, there exists an $x$ such that $L(x, s) \subsetneqq L(x, r)$, and thus there exists a $y \in L(x, r) \setminus L(x, s)$.  Therefore there exists a swing sequence $\seqq{r}$ such that $y \in \beta(x, r_i)$ for all $i$.  Now clearly there cannot exist a $k$ such that $r_k \leq s$, for if there were, taking $s_1 = s$ and $s_i = r_{i+k-1}$ generates a swing sequence for $s$ such that $y \in \beta(x, s_i)$ for all $i$, and this would imply that $y \in L(x, s)$.  Thus $r_k > s$ for all $k$, which completes the proof.} 
\lemma{\label{lemma:levellessthan}If $(X, R, \beta)$ is an ordered $\beta$-space, then if $s <_L r$, 
\begin{enumerate} 
\item $\beta(x, s) \subseteq L(x, r)$ for all $x$ 
\item $L(x, s) \subseteq \beta(x, r)$ for all $x$ 
\item $s< r$ in the ordering on $R$ 
\end{enumerate}} 
\Pf{\begin{enumerate} 
\item Let $\seqq{r}$ be a swing sequence for $r$ such that $ s < r_k$ for all $k$.  Then $\beta(x, s) \subseteq \beta(x, r_k)$ for all $k$, and so $\beta(x, s) \subseteq L(x, r)$ as desired. 
\item This is trivial, since by definition, $L(x, s) \subseteq L(x, r) \subseteq \beta(x, r)$ for all $x$. 
\item Suppose not.  Then $s \geq r$.  Fixing some $x \in X$, for any $y \in L(x, r)$ there is a swing sequence $\seqq{r}$ for $r$ such that $y \in \beta(x, r_i)$ for all $i$.  Now we can obtain a swing sequence $\seqq{s}$ for $s$ by taking $s_1 = s$ and $s_n = r_n$ for $n > 1$.  Clearly $y \in \beta(x, s_i)$ for all $i$, and so $y \in L(x, s)$.  But this implies that for all $x$, $L(x, r) \subseteq L(x, s)$, or that $r \leq_L s$, a contradiction. 
\end{enumerate}} 
\lemma{\label{lemma:levelcontain}Let $(X, R, \beta)$ be an ordered $\beta$-space.  Suppose $\seqq{a} \subseteq L(x, r)$, and let $s \in R$ satisfy $s <_L r$.  Then if $A$ is the $s$-convergent set for $\seqq{a}$, $A \subseteq L(x, r)$.} 
\Pf{Let $y \in A$, and let $\seqq{s}$ be a swing sequence for the convergence to $y$.  Then there exists an $N \in \N$ such that, for all $n \geq N$, $a_n \in \beta(y, s_2)$.  This implies that $\beta(y, s_2) \subseteq \beta(a_N, s_1) $.  But since $s_1 = s <_L r$, by the above lemma we know that $\beta(a_N, s) \subseteq L(a_N, r)$.  Since $a_N \in L(x, r)$, by Proposition \ref{prop:levelequality} we have that \[ y \in \beta(a_N, s) \subseteq L(a_N, r) = L(x, r) \] as desired.} 
  
\Thm{Contraction Mapping Theorem}{Let $(X, R, \beta)$ be Hausdorff, ordered, radially complete, level structured, and CDLB.  Then any contraction $f:X \to X$ has a unique fixed point.} 
\Pf{Pick any $x_1 \in X$ and let $N$ be the contraction degree for $f$.  Let $r_1$ be such that $f(x_1), f^2(x_1), \ldots, f^N(x_1) \in \beta(x_1, r_1)$.  By Proposition \ref{prop:contraction}, the sequence $f^n(x_1)$ is $r_1$-Cauchy, and so since our space is radially complete, the $r_1$-convergent set $A_1$ is nonempty with $f(A_1) \subseteq A_1$.  Inductively suppose that we have $x_k$, $r_k$, and $A_k$ defined, with $f(A_k) \subseteq A_k$.  Pick $x_{k+1} \in A_n$ arbitrarily.  If $f(x_{k+1}) = x_{k+1}$, we set $x^* = x_{k+1}$ and skip down to the section of the proof where we show that $x^*$ is the unique fixed point of $f$.  If $x_{k+1}$ is not a fixed point, then since our space is Hausdorff, $x_{k+1}$ and $f(x_{k+1})$ can be separated by disjoint open sets, and thus, since our space is structured, there exists an $r_{k+1} <_L r_k$ such that \[ f(x_{k+1}), \ldots, f^N(x_{k+1}) \in \beta(x_{k+1}, r_{k+1}) \]  Then as before we define $A_{k+1}$ to be the $r_{k+1}$-convergent set of the sequence $f^n(x_{k+1})$, and again, since our space is radially complete, we know that $A_{k+1}$ is nonempty.  Now by Lemma \ref{lemma:levelcontain}, since $A_{k} = L(x_{k+1}, r_{k})$ and $r_{k+1} <_L r_k$ and $(f^i(x_{k+1}))_{i=1}^{\infty} \subseteq L(x_{k+1}, r_k)$, then we have that \[ A_{k+1} = L(x_{k+2}, r_{k+1}) \subseteq L(x_{k+1}, r_k) \]  Thus, by induction, we have a countable collection \[ A_1 \supseteq A_2 \supseteq \ldots \supseteq A_n \supseteq \ldots \] with the properties established above.  We now claim that the sequence $\seqq{x}$ is Cauchy.  To see this, pick any $s \in R$, and let $t$ be a swing value for $s$.  Now, since our space is CDLB and the sequence $\seqq{r}$ satisfies $r_{i+1} <_L r_i$ for all $i$, then the collection $\seqq{r}$ is a countable level-base, and so there must be some $r_k <_L t$.  But for all $n, m \geq k$, \[ x_{n+1} \in A_k = L(x_{k+1}, r_k) \subseteq \beta(x_{k+1}, t) \]  Thus, for all $n, m \geq k+1$, \[ x_{n+1} \in \beta(x_k, t) \subseteq \beta(x_{m+1}, s) \] as desired.\\~~\\ 
\par 
Since our space is radially complete, it is complete, and so our Cauchy sequence $\seqq{x}$ must converge to a point $x^*$.  Now we show that $f(x^*) = x^*$.  Suppose not.  Then since our space is Hausdorff, there is an $r \in R$ such that $f(x^*) \notin \beta(x^*, r)$.  Let $s$ be a swing value of $r$ and let $t$ be a swing value of $s$.  Now choose $N$ such that $x^*, f^n(x_N) \in \beta(x_{N+1}, t)$ for all $n$.  Then \[ x^* \in \beta\left( f^n(x_N), s \right) \] for all $n$.  But since $f$ is a contraction, this gives us that \[ f(x^*) \in f \left( \beta \left( f^n(x_N), s \right) \right) \subseteq \beta \left( f^{n+1}(x_N), s \right) \subseteq \beta(x^*, r) \] which is a contradiction.\\~~\\ 
\par 
Finally, we show that this fixed point is unique.  Suppose that $y^* \neq x^*$ is another fixed point of $f$.  Let \[ G = \left\{ r \in R \; : \; y^* \in \beta(x^*, r) \right\} / =_L \]  Since our space is structured, $G$ cannot be empty.  Then $G$ must contain a minimal element, since if not, we can find a countable strictly-decreasing collection $([r_i])_{i=1}^{\infty} \subseteq G$ which, since our space is CDLB, must be a countable level-base.  But since our space is Hausdorff there must be some radial value $t$ where $y^* \notin L(x^*, t)$, so this is a contradiction.  Now let $[k] \in G$ be the minimal element, and let the representative $k$ be such that $y^* \in \beta(x^*, k)$.  Then as before, by taking $f^i\left( \beta(x^*, k) \right)$, we know that these iterates generate a swing sequence $\seqq{k}$ of $k$ where, since $x^*$ is a fixed point, \[ f^{iN}\left( \beta(x^*, k) \right) \subseteq \beta(x^*, k_i) \]  But since $y^*$ is also a fixed point and $y^* \in \beta(x^*, k)$, we have that $ y^* \in f^i\left( \beta(x^*, k) \right)$ for all $i$, and thus that \[ y^* \in \beta(x^*, k_i) \] for all $i$.  But then $y^* \in L(x^*, k)$, and so by structuredness, there is a $k' <_L k$ such that $y^* \in \beta(x^*, k')$.  But this contradicts that $[k]$ was the minimal element in $G$, and so such a $y^*$ cannot exist.} 
\prop{Any field-metric space is structured, and if it satisfies that $s$ is a swing value of $r$ implies that $s \leq r/2$, then it is level structured.} 
\Pf{Given $x \in X$ and $(y_i)_{i=1}^n \subseteq X$, we set $m = \max\{ d(x, y_i) \; : \; y = 1, 2, \ldots, n \}$ where $d$ is the field-metric on $X$.  If $m > 0$, $(y_i)_{i=1}^n \subseteq \beta(x, 2m)$.  If $m=0$, then $y_i = x$ for all $i$, so that $(y_i)_{i=1}^n \subseteq \beta(x, 1)$, and so our space is structured. \\~~\\ 
\par 
Now suppose that $(y_i)_{i=1}^n \subseteq L(x, r)$ for some $r \in R$.  Further suppose that there is some $y_k$ such that $y_k$ can be separated by $x$ via disjoint open neighborhoods, and note that this condition is equivalent to the statement that $m > 0$.  Now for any swing sequence $\seqq{r}$ of $2m$, necessarily $r_2 \leq m$ by our hypothesis.  Thus, if $k$ is such that $d(x, y_k) = m$, since $m > 0$, $y_k \notin \beta(x, m)$, and thus $y_k \notin \beta(x, r_2)$.  This implies that $y_k \notin L(x, 2m)$, and since our ordered space is necessarily level-ordered, we have that $2m <_L r$ as desired.} 
\thm{The field-metric $\L$ of formal Laurent series over a radially complete linearly ordered field is radially complete.} 
\Pf{REQUIRES PROOF.  But the concept is simple.  When we showed completeness, we showed that a Cauchy sequence had the ``slot machine'' convergence.  The same is true with $r$-convergence, except that the slot machine stops at position $n$, where $n$ is such that $r = \sum_{i=n}^{\infty} a_i x^i$.  Can probably just move / copy the proof from section 1 here and amend it.} 
\prop{The field-metric $\L$ of formal Laurent series over a linearly ordered field is level structured.} 
\Pf{To show that $\L$ is level structured, it suffices to show that for any $r \in R$, $s$ is a swing value of $r$ implies that $s \leq r/2$.  Suppose that $s > r/2$.  Then for any $x$ there exists a $y \in \beta(x, s)$ such that $d(x, y) = |x-y| = (2s+r)/4$.  Let $z = x - y$, and set $y' = x + z$.  Then $d(x, y') = |x - y'| = |z| = (2s+r)/4$ so that $y' \in \beta(x, s)$ as well.  But \[ d(y, y') = |y - y'| = |(x+z) - (x-z)| = |2z| = (2s+r)/2 = s + r/2 > r \] so that $y' \notin \beta(y, r)$.  In general this means that $\beta(x, s) \not\subseteq \beta(y, r)$, so that $s$ is not a swing value of $r$.} 
\prop{The field-metric $\L$ of formal Laurent series over an archimedean linearly ordered field is CDLB.} 
\Pf{First we show that for any $r = \sum_{i=n}^{\infty} a_i x^i$, $L(y, r) = L(y, x^n)$ for any $y$.  From part (2) above, if $\seqq{r}$ is a swing sequence for $r$, then $r_i \leq r \cdot 2^{1-i}$.  Therefore, \[ L(y, r) = \{ z \; : \; z \in \beta(y, r/2^i)  \textrm{ for all } i \in \N \} \]  Now we will show that $L(y, r) \subseteq L(y, x^n)$.  Let $z \in L(y, r)$.  Then $z \in \beta(y, r/2^i)$ for all $i$.  But since the base field for $\L$ is archimedean, then there is an $N$ such that $a_n / 2^N < 1$, so that $r / 2^N < x^n$.  Therefore the sequence $\seqq{\alpha}$, where $\alpha_1 = x^n$ and $\alpha_i = r/2^{N+i}$ when $i > 1$ is a swing sequence for $x^n$, and therefore $z \in L(y, x^n)$ as desired.  The reverse argument suffices to show the opposite inclusion, where we simply note that since our field is archimedean, there is an $N$ such that $2^{-N} < a_n$.} 
\cor{Let $\L$ denote the field of formal Laurent series over $\R$.  Then any function $f:\L \to \L$ such that $d\left( f(y), f(z) \right) \leq r \cdot d(y, z)$ for some $r$ infinitesimally close to an element of $[0, 1) \subseteq \R$ has a unique fixed point.} 
\Pf{We note that $\R$ is archimedean, linearly ordered, and radially complete.  Therefore, $\L$ is CDLB, level structured, and radially complete.  Since $\L$ is a field metric, it is Hausdorff and ordered.  Therefore to complete the proof it suffices to show that a function $f$ satisfying the given condition is necessarily a contraction.  First we note that since our $r$ is infinitesimally close to an element of $[0, 1) \subseteq \R$, necessarily $r \in [0, 1) \subseteq \L$, so that \[ f\left( \beta(y, s) \right) \subseteq \beta\left( f(y), s \right) \] trivially.  Now since $r$ is infinitesimally close to some value in $[0, 1) \subseteq \R$, let $r'$ be that value, and let $s = (1+r')/2$.  We note that $r < s < 1$ in the ordering on $\L$, and further that $s \in (0, 1) \subseteq \R$.  Therefore, there is an $N$ such that $s^N \leq 1/2$, which implies that $r^N \leq 1/2$.  Since we know that \[ d\left( f^N(y), f^N(z) \right) \leq r^N \cdot d(y, z) \] then it follows that \[ f^N \left( \beta(y, t) \right) \subseteq \beta \left(f^N(y), t/2 \right) \] so that $f$ is a contraction as desired.} 
  
\prop{The field-metric space $\L^n$ of formal Laurent series is Hausdorff, ordered, radially complete, level-structured, and CDLB, where we take as field-metric the function $d((a_i)_{i=1}^n,(b_i)_{i=1}^n)=\sup_i\{|b_i-a_i|\}$.} 
\Pf{Since $\L^n$ is a field-metric space, it is Hausdorff and ordered. To show it is radially complete, fix $r\in R=(\L^n)^+$. Let $(a_i)_{i=1}^\infty$ be an $r$-Cauchy sequence in $\L^n$, where $a_i=(\sum_{k=M}^\infty b_{i,j,k}x^k)_{j=1}^n$ for each $i$ (note that since the sequence is $r$-Cauchy, there is a finite lower bound $M$ on the set of all indices). Since the sequence is $r$-Cauchy there exists a swing sequence $(r_i)_{i=1}^n$ such that for all $i$ there exists an $N_i\in\N$ such that if $l,m\geq N_i$, then $d(a_l,a_m)<r_i$; that is, \[\left|\left(\sum_{k=M}^\infty b_{l,j,k}x^k\right)-\left(\sum_{k=M}^\infty b_{m,j,k}x^k\right)\right|<r_i\] for all $j$ and for all $l,m\geq N_i$. Let $P=\{p\in\Z\ |\ \exists\ r_i\leq x^p\}$. 
First consider the case where $P$ is bounded above, i.e. there exists a $q\in\Z$ such that $p<q$ for all $p\in P$, and pick the smallest such $q$. For $k=M,\ldots,q-1$, pick $i_k$ to be the smallest integer such that $r_{i_k}\leq x^k$. Then if $l,m\geq N_{i_k}$, we have \[\left|\left(\sum_{k=M}^{\infty} b_{l,j,k} x^k\right)-\left(\sum_{k=M}^{\infty} b_{m,j,k} x^k\right)\right|<r_{i_k}\leq x^k\] for all $j$. Then $b_{l,j,k}=b_{m,j,k}$ for all $l,m\geq N_{i_K}$, for all $j$, and for all $k\leq K-1$. Let $d_{j,k}=b_{N_{i_{k+1}},k}$ for all $j$ and for all $k\leq q-2$. 
For all $i\geq i_{q-1}$ write $r_i=\sum_{k=q-1}^\infty \alpha_{i,k}x^k$ for some $\alpha_{i,k}\in\R$. Since $(r_i)$ is a swing sequence, we must have $\alpha_{i+1,q-1}\leq\alpha_{i,q-1}/2$ for all $i\geq i_{q-1}$, so $\alpha_{i,q-1}\to 0$ in $\R$ as $i\to\infty$. Thus there exists a $d_{j,q-1}\in\R$ such that $b_{i,j,q-1}\to d_{j,q-1}$ in $\R$ as $i\to\infty$. For $k\geq q$ and for all $j$, let $d_{j,q}\in\R$. Let $c=(\sum_{k=M}^\infty d_{j,k}x^k)_{j=1}^n$. Then $a_i\stackrel{r}{\to}c$. 
Now consider the case where $P$ is unbounded, i.e.\ there does not exist a $q$ such that $p<q$ for all $p\in P$. In this case, the sequence $(a_i)$ is actually Cauchy. Since $\L$ is complete, there exists a $c\in\L^n$ such that $a_i\to c$, and thus $a_i\stackrel{r}{\to}c$. Therefore $\L^n$ is radially complete. 
To show $\L^n$ is level structured, let $z,y_1,\ldots,y_m\in\L^n$, where $z=(\sum_{k=M}^\infty a_{j,k}x^k)_{j=1}^n$ and $y_i=(\sum_{k=M}^\infty b_{i,j,k}x^k)_{j=1}^n$ for each $i$, with at least one of the $a_{j,M}$ and $b_{i,j,M}$ is nonzero. Then $y_i\in L(z,x^{M-1})$ for all $i$. 
Let $\alpha=\max_{i,j}\{|a_{j,M}-b_{i,j,M}|\}$. Then $2\alpha x^M<_L x^{M-1}$, and $y_i\in\beta(z,2\alpha x^M)$ for all $i$. Thus $\L^n$ is level structured. 
Finally, to show $\L^n$ is CDLB, let $(r_i)_{i=1}^\infty$ be a sequence in $R$ such that $r_{i+1}<_L r_i$ for all $i$. Let $s\in R$. Write $r_i=\sum_{k=M}^\infty \alpha_{i,k}x^k$ and $s=\sum_{k=N}^\infty \beta_k x^k$, where $\alpha_{1,M},\beta_N\neq 0$. Let $m=\min\{N-M+1,1\}$. Then $r_m<_L s$, and thus $\L^n$ is CDLB.}

\end{document}